\documentclass[10pt]{amsart}

\usepackage{amssymb}

\usepackage{enumerate}

\usepackage{graphicx}

\makeatletter
\@namedef{subjclassname@2010}{%
  \textup{2010} Mathematics Subject Classification}
\makeatother

\newtheorem{thm}{Theorem}[section]
\newtheorem{corollary}[thm]{Corollary}
\newtheorem{lemma}[thm]{Lemma}
\newtheorem{proposition}[thm]{Proposition}

\theoremstyle{definition}

\newtheorem{remark}[thm]{Remark}
\newtheorem{example}[thm]{Example}

\numberwithin{equation}{section}

\begin{document}


\baselineskip=17pt



\title[3D Corner manifolds]{Corner manifolds}

\author[Beldjilali Gherici]{Beldjilali Gherici}
\address{ Laboratory of Quantum Physics and Mathematical Modeling (LPQ3M)\\
    University of Mascara,   Algeria.}
\email{gherici.beldjilali@univ-mascara.dz}

\date{}

\begin{abstract}
The aim of this paper is two-fold. First, the study of  $C_{12}$-structure (called by us corner structure) is extended to the general case without any condition, unlike our previous papers (see, \cite{BB, BG2, BG, BBB}). Second, starting from a $C_{12}$-structure which is not normal, we construct Trans-Sasakian structures which are normal. Class of examples are given.
\end{abstract}



\maketitle

\section{Introduction}

In \cite{CG}, D. Chinea and C. Gonzalez have obtained a classification of the almost contact metric structures, studying the space that possess the same symmetries as the covariant derivative of the fundamental 2-form. This space is decomposed into twelve irreducible components $C_1,...,C_{12}$.\\
In dimension 3, the classes $C_i$ reduce to the classes: $\vert C \vert$ class of cosymplectic structures, $C_5$ class of $\beta$-Kenmotsu structures, $C_6$ class of $\alpha$-Sasakian structures, $C_9$-structures that are almost cosymplectic structures and $C_{12}$-structures which are the structures targeted in this work.

Most of the research related to almost contact metric structures is concerned with the normal structures which contain the first three classes. For the $C_{12}$ class which is not normal, some studies have recently started to appear on this subject,  see for-example
\cite{ BB, BG2, BG, BBB, CF}.

The present paper is devoted to 3-dimensional $C_{12}$-manifolds. We present a detailed study of such manifold  in the third dimension and highlight some of its geometric properties with concrete examples. We will also explain the reason for calling these manifolds " Corner manifolds". Next, we give a direct method to pass from the corner manifold to cosymplectic and Kenmotsu manifold. More generally, we investigate a new deformation  of the structural tensor $\varphi$ and metric tensor $g$ at the same time, which allow us to construct a Trans-Sasakian manifold starting from a corner structure. 

First of all, we will start by introducing the basic concepts that we need in this research.

\section{Preliminaries}
An odd-dimensional Riemannian manifold $(M^{2n+1},g)$ is said to be an almost contact metric manifold if there exist on $M$ a $(1,1)$-tensor field $\varphi$, a vector field $\xi$ (called the structure vector field) and a $1$-form $\eta$ such that 
\begin{equation}\label{CondPresqCont}
 \left\{
          \begin{array}{lll}
\eta(\xi)=1,\\
\varphi^{2}(X) = -X+\eta(X)\xi ,\\
g(\varphi X,\varphi Y) = g(X,Y)-\eta(X)\eta(Y),
          \end{array}
  \right.
\end{equation}
 for any vector fields $X$, $Y$ on $M$. 
 
 In particular, in an almost contact metric manifold we also have
 $$ \varphi\xi=0 \qquad  and \qquad \eta \circ \varphi=0.$$

The fundamental 2-form $\phi$ is defined by 
$$\phi (X,Y) = g(X, \varphi Y).$$
It is known that the almost contact
structure $(\varphi, \xi, \eta)$  is said to be normal if and only if
\begin{equation}\label{N1phi}
N^{(1)}(X,Y)=N_{\varphi}(X,Y)+2{\rm d}\eta(X,Y)\xi=0,
\end{equation}
 for any $X$, $Y$ on $M$, where where ${\rm d}$ denotes the exterior derivative and $N_{\varphi}$  the Nijenhuis torsion of $\varphi$, given by
\begin{equation}\label{Nphi}
N_{\varphi}(X,Y)=\varphi^{2}[X,Y]+[\varphi X,\varphi Y]-\varphi[\varphi X,Y]-\varphi[X,\varphi Y].
\end{equation}

In \cite{OLS}, Olszak completely characterized the local nature of the normal almost contact metric manifold of dimension three.  He  proved that, for an arbitrary $3$-dimensional almost contact metric manifold $(M, \varphi, \xi, \eta,g)$, we have
\begin{equation}\label{PropOls}
 \left\{
          \begin{array}{lll}
&(1)& (\nabla_X \varphi )Y=g(\varphi \nabla_X \xi , Y)\xi - \eta(Y)\varphi \nabla_X \xi\\
&(2)&{\rm d}  \phi = 2 \beta\eta \wedge \phi \\
&(3)& {\rm d} \eta = \eta \wedge ( \nabla_{\xi} \eta ) +\alpha  \phi,
          \end{array}
  \right.
\end{equation}
where $2\alpha =\big(tr_g(\varphi \nabla \xi )\big)$, $2\beta={ \rm div} \xi $
and  $\nabla$ denotes the Levi-Civita connection of $g$.

Also, He show that any   almost contact metric manifold of dimension three is normal if and only if
\begin{equation}\label{DifTrans-Sas}
 \nabla_X  \xi = -\alpha \varphi X - \beta \varphi^2 X.
\end{equation}
This mean that any normal almost contact metric manifold of dimension three  is a Trans-Sasakian manifold of type ($\alpha$ , $\beta$)  introduced by  J. A. Oubi$\tilde{n}$a \cite{OUB}.

\section{Three dimensional  $C_{12}$-manifold} 

In the classification of D. Chinea and C. Gonzalez \cite{CG} of almost contact metric manifolds there is a class called $C_{12}$-manifolds  which can be integrable but never normal. In this classification, $C_{12}$-manifolds are defined by
\begin{equation}\label{CGC12}
(\nabla_X \phi)(Y,Z)=\eta(X)\eta(Z) (\nabla_{\xi} \eta) \varphi Y - \eta(X)\eta(Y) (\nabla_{\xi} \eta) \varphi Z.
\end{equation}

In \cite{BBB} and \cite{CF}, The $(2n+1)$-dimensional $C_{12}$-manifolds are characterized by:

\begin{equation}\label{NablaVarphi}
 (\nabla_X \varphi) Y =  \eta(X) \big( \omega(\varphi Y) \xi + \eta(Y) \varphi \psi \big),
\end{equation}
for any $X$ and $Y$ vector fields on $M$, where  $ \omega = -\big( \nabla_{\xi}\xi \big)^{\flat} = -\nabla_{\xi}\eta$ and $\psi$ is the vector field given by $$\omega(X)=g(X, \psi) = -g(X, \nabla_{\xi}\xi).$$

Moreover,  in  \cite{BBB} the $(2n+1)$-dimensional $C_{12}$-manifolds is also  characterized by 
\begin{equation}\label{cond11}
 {\rm d} \eta =   \omega \wedge \eta  \qquad  {\rm d} \phi = 0 \qquad  and \qquad N_{\varphi} =0.
\end{equation}
\begin{remark}
Unlike our previous papers related to $C_{12}$-manifold, the existence of three global vector fields $\xi$, $\psi$ and $\varphi \psi$ which are naturally pairwise orthogonal, encourages us to suggest the name "\textbf{Corner manifold}" instead of $C_{12}$-manifold.
\end{remark}

In dimension three, two nice characterizations are given in \cite{BG2}:
\begin{thm}\label{ThNablaXi}
Let $(M^3, \varphi, \xi, \eta,g)$ be a $3$-dimensional almost contact metric manifold. $M$ is a corner manifold if and only if 
\begin{equation}\label{NablaXi}
\nabla_X \xi = - \eta(X)\psi.
\end{equation}
Or, equivalently
\begin{equation}
\nabla_{\varphi X} \xi = 0.
\end{equation}
\end{thm}
All of the above definitions of the $C_{12}$-Manifolds (\ref{CGC12}), (\ref{NablaVarphi}) and  (\ref{cond11}) have been shown to be equivalent (see \cite{BG2, BBB, CF}).

In  \cite{BB, BG2, BG, BBB}, the authors studied the 3-dimensional unit corner manifold with $\omega$ is a closed 1-form on $M$  i.e. the case where $\psi$ is a unit vector field and ${\rm d}\omega=0$. 

More precisely, in \cite{BBB} the authors studied the 3-dimensional  $C_{12}$-structures with $\vert \psi \vert=1$ and ${\rm d}\omega=0$ and  in \cite{BB, BG2,  BG} the  3-dimensional $C_{12}$-structure was investigated under the condition ${\rm d}\omega=0$ only.   Here, we will deal with the general case, i.e. $\psi$ and $\omega$ are not necessarily unitary and closed respectively. For that, taking $V = {\rm e}^{-\rho}\psi$ where $ {\rm e}^{\rho} = \vert \psi \vert$, we get immediately that  $\{\xi, V, \varphi V \}$ is  an orthonormal frame. We refer to this basis as  fundamental basis. Consequently, $\{\eta , \theta^1 = V^{\flat}={\rm e}^{-\rho} \omega , \theta^2 = - {\rm e}^{-\rho} \omega \circ \varphi\}$
is  the fundamental dual frame we will use it later. 

Now, using the fundamental basis, we can get the following:
\begin{proposition}\label{PropِConnex}
For any $C_{12}$-manifold, for all vector field $X$ on $M$ we have
\begin{itemize}
  \item[1)] $\nabla_{X} \xi =  - {\rm e}^{\rho}\eta(X) V,$
  \item[2)] $\nabla_{\xi }V= {\rm e}^{\rho} \xi +  \sigma  \varphi V,
                                         \qquad\qquad where \quad \sigma =g(\nabla_{\xi }V, \varphi V)$
  \item[3)] $ \nabla_V  V = \varphi V(\rho) \varphi V,$
  \item[4)] $ \nabla_{\varphi V}V =\big( {\rm div} V -{\rm e}^{\rho}\big) \varphi V,$
  \item[5)] $ \nabla_\xi \varphi V = -\sigma V,$
  \item[6)] $ \nabla_V  \varphi V =-\varphi V(\rho) V, $ 
  \item[7)] $ \nabla_{\varphi V}\varphi V =\big( {\rm e}^{\rho} - {\rm div} V \big) V.$
 \end{itemize}
\end{proposition}
\begin{proof}
For the first, just use (\ref{NablaXi}) with  $\psi = {\rm e}^{\rho} V$. 
 For the second, we have 
\begin{align*}
\nabla_{\xi }V &= g(\nabla_{\xi }V, \xi)\xi + g(\nabla_{\xi }V, V)V + g(\nabla_{\xi }V, \varphi V)\varphi V\\
&= - g(V, \nabla_{\xi }\xi)\xi + g(\nabla_{\xi }V, \varphi V)\varphi V\\
&={\rm e}^{\rho}\xi + \sigma  \varphi V,
\end{align*}
where $\sigma =g(\nabla_{\xi }V, \varphi V)$. 

For $\nabla_{V} V$, we have $ {\rm d}\eta= \omega \wedge \eta$ implies ${\rm d} \omega \wedge \eta =0$. Then
\begin{align*}
0 &=  ({\rm d} \omega \wedge \eta)(X, \psi , \xi) \\
&= {\rm d} \omega(X, \psi) + \eta(X) {\rm d} \omega(\psi, \xi)\\
&= g(\nabla_X \psi , \psi) - g(\nabla_{\psi}\psi , X) + \eta(X)\big( g(\nabla_{\psi}\psi , \xi) - g(\nabla_{\xi}\psi , \psi)  \big)\\
&= {\rm e}^{2\rho} X(\rho) -g(\nabla_{\psi}\psi , X)+ \eta(X)\big( g(\psi , \nabla_{\psi}\xi) - {\rm e}^{2\rho} \xi(\rho) \big).
\end{align*}
Therefore
\begin{align*}
g(\nabla_{\psi}\psi , X) = {\rm e}^{2\rho} X(\rho) - {\rm e}^{2\rho} \xi(\rho)\eta(X),
\end{align*}
i.e. $\nabla_{\psi} \psi ={\rm e}^{2\rho} {\rm grad}\rho  - {\rm e}^{2\rho} \xi(\rho)\xi$ which gives $ \nabla_V V = {\rm grad}  \rho - \xi(\rho)\xi -V(\rho) V,$\\
since
\begin{align*}
{\rm grad}  \rho &= \xi(\rho) \xi + V(\rho) V + \varphi V(\rho) \varphi V,
\end{align*}
then,
$$ \nabla_V V = \varphi V(\rho) \varphi V.$$
For the forth, we have
\begin{align*}
\nabla_{\varphi V} V  &= g(\nabla_{\varphi V} V , \xi) \xi + g( \nabla_{\varphi V} V ,V) V + g(\nabla_{\varphi V} V , \varphi V) \varphi V\\
  &=  g(\nabla_{\varphi V} V , \varphi V) \varphi V.
\end{align*}
Knowing that
\begin{align*}
{\rm div} V  &= g(\nabla_{\xi} V , \xi) + g( \nabla_{V} V ,V) + g(\nabla_{\varphi V} V , \varphi V)\\
  &=  {\rm e}^{\rho} + g(\nabla_{\varphi V} V , \varphi V),
\end{align*}
then $ g(\nabla_{\varphi V} V , \varphi V) = {\rm div} V -{\rm e}^{\rho}$ which gives
$$ \nabla_{\varphi V} V = \big( {\rm div} V -{\rm e}^{\rho}\big) \varphi V.$$
For the rest, just use the formula $\nabla_X \varphi Y =( \nabla_X \varphi )Y + \varphi \nabla_X Y$ with the help of (\ref{NablaXi}). 
\end{proof}

\begin{proposition}\label{Propِd1Forms}
Let $(M, \varphi, \xi,  \eta, g)$ be a three dimensional corner manifold  and $\phi$ its fundamental 2-form. Then,  we have
\begin{itemize}
 \item[1)] $\phi = 2  \theta^2\wedge \theta^1.$
  \item[2)] $ {\rm d} \theta^1 =  \sigma  \eta \wedge \theta^2 + \varphi V(\rho) \theta^1 \wedge \theta^2. $
  \item[3)] ${\rm d} \theta^2=  -\sigma  \eta \wedge \theta^1 -({\rm e}^{\rho} -  {\rm div}V) \theta^1 \wedge \theta^2. $
 \end{itemize}
\end{proposition}
\begin{proof}
Since $\{\eta , \theta^1 = V^{\flat}={\rm e}^{-\rho} \omega , \theta^2 = - {\rm e}^{-\rho} \omega \circ \varphi\}$
is  the dual frame of differential 1-forms on a 3-dimensional corner manifold then, $\{\eta \wedge \theta^1, \eta \wedge \theta^2, \theta^1 \wedge \theta^2\}$
is  the dual frame of differential 2-forms on $M$. \\
$1)$ For all $X$ and $Y$ vector fiels on $M$, we have
\begin{eqnarray*}
\varphi Y &=& g(\varphi Y , \xi)\xi + g(\varphi Y , V)V+ g(\varphi Y , \varphi V)\varphi V \\
&=&- \theta^2(Y)V + \theta^1(Y)\varphi V.
\end{eqnarray*}
Thus,
\begin{eqnarray*}
\phi(X,Y)&=& g(X, \varphi Y)\\
&=& - \theta^2(Y)\theta^1(X) + \theta^1(Y)\theta^2(X)\\
&=& 2  (\theta^2\wedge \theta^1)(X,Y).
\end{eqnarray*}
$2)$ We can write
$$ {\rm d} \theta^1  = a \eta  \wedge \theta^1 + b \eta \wedge \theta^2 + c \theta^1 \wedge \theta^2,$$
where $a, b, c$ are three functions on $M$. Using Proposition \ref{PropِConnex}, we get
\begin{align*}
a  &= 2{\rm d} \theta^1(\xi, V) \\
&= g(\nabla_{\xi}V , V) - g(\nabla_V V , \xi)\\
&=  0,
\end{align*}
and 
$$ b  = 2{\rm d} \theta^1(\xi, \varphi V)= \sigma  ,\qquad c  = 2{\rm d} \theta^1(V, \varphi V)=\varphi V(\rho) .$$
Then,
\begin{align}\label{dOmega}
 {\rm d} \theta^1 &=  \sigma  \eta \wedge \theta^2 + \varphi V(\rho) \theta^1 \wedge \theta^2.
 \end{align}
$3)$ With the same techniques, we obtain the third formula.
\end{proof}
\section{Class of examples}
We denote the Cartesian coordinates in a $3$-dimensional Euclidean space $\mathbb{R}^3$ by $(x_1, x_2, x_3)$ and define a symmetric tensor field $g$ by
\begin{eqnarray*}
g= \left(
       \begin{array}{ccc}
        \tau(x_1,x_2,x_3)^2& 0 & 0\\
        0 &  \kappa(x_1,x_2,x_3)^2 & 0\\
       0& 0 & \mu(x_1,x_2,x_3)^2
       \end{array}
\right),        
\end{eqnarray*}
where $ \tau$, $\kappa$ and $ \mu$ are  functions on $\mathbb{R}^3$  and $ \tau \kappa \mu \neq 0$ everywhere.
Further, we define an almost contact structure  $(\varphi ,\xi ,\eta)$ on $\mathbb{R}^3$ by
$$
\varphi= \left(
       \begin{array}{ccc}
        0 & 0 & 0\\
        0 &  0 &-\frac{\mu}{\kappa}\\
       0 & \frac{\kappa}{\mu}& 0
       \end{array}
\right),\qquad
 \xi= \frac{1}{\tau}\left(
       \begin{array}{lll}
        1\\
        0\\
        0
       \end{array}\right) ,
        \qquad \eta =\big(\tau , 0 , 0 \big).$$
Notice that
${\rm d}\eta = \tau_2 d x_2 \wedge dx_1 + \tau_3 d x_3 \wedge dx_1 = \omega \wedge \eta$ with 
$ \omega = \frac{\tau_2}{\tau}dx_2 + \frac{\tau_3}{\tau} dx_3,$
where
$ \tau_i = \frac{\partial \tau}{\partial x_i}$ (The index always indicates the derivative). 
Therefore
$$\psi = \frac{\tau_2}{\tau \kappa^2 }\frac{\partial}{\partial x_2}  + \frac{\tau_3}{\tau \mu^2 }\frac{\partial}{\partial x_3}.$$ 
Now, we can see that 
$$\vert \psi \vert ^2=  \frac{1}{\tau^2}\Big( \frac{\tau_2^2}{\kappa^2} + \frac{\tau_3^2}{\mu^2}\Big)\qquad and \qquad
{\rm d}\omega= \Big(\frac{\tau_2}{\tau}\Big)_1 dx_1 \wedge dx_2 + \Big(\frac{\tau_3}{\tau}\Big)_1 dx_1 \wedge dx_3.$$
We give the following orthonormal basis
$$ \xi = \frac{1}{\tau}\frac{\partial}{\partial x_1} ,\qquad e_1= \frac{1}{\kappa}\frac{\partial}{\partial x_2},\qquad e_2=\frac{1}{\mu} \frac{\partial}{\partial x_3}.$$
So, the components of the Levi-Civita connection corresponding to $g$ are written
$$
      \begin{array}{ccc}
        \nabla_{ \xi} \xi= -\frac{\tau_2}{ \kappa \tau}e_1 - \frac{ \tau_3}{\tau \mu}e_2,
         & \nabla_{ \xi}e_1 =\frac{\tau_2}{\kappa \tau } \xi,
         &  \nabla_{ \xi}e_2= \frac{\tau_3}{ \tau \mu} \xi, \\
        \nabla_{e_1} \xi =\frac{\kappa_1}{\kappa \tau}e_1,
        &  \nabla_{e_1}e_1=- \frac{ \kappa_1}{\kappa \tau} \xi  -\frac{\kappa_3}{ \kappa \mu}e_2,
        & \nabla_{e_1}e_2=\frac{\kappa_3}{ \kappa \mu} e_1,\\
         \nabla_{e_2} \xi =\frac{\mu_1}{\tau \mu}e_2,
        &  \nabla_{e_2}e_1= \frac{\mu_2}{\kappa \mu}e_2,
        & \nabla_{e_2}e_2=-\frac{\mu_1}{ \tau \mu}  \xi-\frac{\mu_2}{ \kappa  \mu} e_1.
        \end{array}
$$
Using Theorem \ref{ThNablaXi}, one can check that $(\mathbb{R}^3, \varphi ,\xi ,\eta, g)$ is a $3$-parameter family of $C_{12}$-manifolds if and only if 

$$\nabla_{e_i} \xi =- \eta(e_i)\psi =- \eta(e_i)\Big(  \frac{\tau_2}{\kappa \tau}e_1+\frac{\tau_3}{\mu\tau} e_2 \Big),$$
    where   $i \in \{0,1,2\}$    with
     $e_0=\xi$,
    i.e.
    $$ \nabla_{ \xi} \xi= -\frac{\tau_2}{ \kappa \tau}e_1 - \frac{ \tau_3}{\tau \mu}e_2,\qquad  \nabla_{ e_1} \xi=  \nabla_{ e_2}\xi =0.$$
       From the above  components of the Levi-Civita connection, we get
      $$ \kappa_1 = \mu_1 =0.$$

\section{From 3D corner structure to others}
According to the previous section, here we present the sufficient and necessary techniques to construct  normal almost contact metric structures,  especially Sasakian, Kenmodsu and cosymplectic manifolds starting from a corner manifold which is non normal. 
\subsection{The twin structures of corner structure}
On an arbitrary oriented Riemannian 3-manifold one can canonically define a cross product $\times$ of two vector fields
$ X$ and $Y$ on $M$ as follows:
\begin{equation}
g(X \times Y, Z) = dv_g(X, Y, Z),
\end{equation}
for any vector fields $Z$ on $M$. where $dv_g$  denotes the volume form defined by $g$. When $(M, \varphi ,\xi ,\eta, g)$ is an almost contact metric 3-manifold, the cross product is given by the formula \cite{CFG, CAM}:
\begin{equation}
X \times Y= g( \varphi X, Y)\xi - \eta(Y)\varphi X + \eta(X)\varphi Y,
\end{equation}
for all $ X$ and $Y$ vector fields on $M$. Easily, we can notice that $\varphi X = \xi \times X$. 

For our first application, let's define a $(1,1)$-tensor field $\overline{\varphi}$ on $M$ by
$$ \overline{\varphi}X=V \wedge X ,$$
that is, for all $X$ and $Y$  vector fields on $M$
\begin{eqnarray*}
\overline{\varphi}X &=&g( \varphi V, X)\xi - \eta(X)\varphi V\\
&=&\theta^2(X)\xi- \eta(X)\varphi V.
\end{eqnarray*}
\begin{proposition}\label{ACMS1}
Let $(M, \varphi, \xi,  \eta, g)$ be a three dimensional  corner manifold. Then,  $(M,\overline{\varphi}, V,  \theta^1, g)$
is an almost contact metric manifold.
\end{proposition}
\begin{proof}
Easily we can see that $ \overline{\varphi} V =0$ and for all vector field $X$ on $M$, we have
\begin{align*}
\overline{\varphi}^2 X &=\theta^2(X)\overline{\varphi}\xi- \eta(X)\overline{\varphi}\varphi V\\
&= -\theta^2(X)\varphi V - \eta(X)\xi\\
&= -X + \theta^1(X)V,
\end{align*}
because
\begin{align*}
X &=g(X, \xi)\xi+ g(X, V)V + g(X, \varphi V)\varphi V\\
&= \eta(X)\xi + \theta^1(X)V +\theta^2(X)\varphi V.
\end{align*}
For the condition of compatibility, we have
\begin{align*}
g( \overline{\varphi}  X , \overline{\varphi}  Y)  &= \theta^2(X)\theta^2(Y)+ \eta(X)\eta(Y)\\
&= g \big( \theta^2(X)\varphi V+ \eta(X)\xi ,Y \big)\\
&=g \big( X - \theta^1(X)V ,Y \big)\\
&= g(X,Y) - \theta^1(X)\theta^1(Y).
\end{align*}
which completes the proof. 
\end{proof}
Based on these facts, we give the following Theorem:
\begin{thm}\label{ThKen}
The manifold $(M,\overline{\varphi}, V,  \theta^1, g)$ defined above is a $\beta$-Kenmotsu where $\beta =  {\rm e}^{\rho}$
if and only if ${\rm div} V =2{\rm e}^{\rho}$ and $\sigma=\varphi V(\rho)=0$.
\end{thm}
\begin{proof}
From Proposition \ref{PropِConnex}, we have
\begin{equation}\label{SysKen1}
 \left\{
          \begin{array}{lll}
\nabla_{\xi }V= {\rm e}^{\rho} \xi +  \sigma  \varphi V,\\
\nabla_V  V = \varphi V(\rho) \varphi V,\\
\nabla_{\varphi V}V =\big( {\rm div} V -{\rm e}^{\rho}\big) \varphi V.
          \end{array}
  \right.
\end{equation}
On the other side, we have
\begin{equation}\label{SysKen2}
 \left\{
          \begin{array}{lll}
 -\alpha \overline{\varphi} \xi - \beta \overline{\varphi}^2\xi= \beta \xi +  \alpha \varphi V,\\
 -\alpha \overline{\varphi} V - \beta \overline{\varphi}^2V = 0,\\
 -\alpha \overline{\varphi} \varphi V - \beta \overline{\varphi}^2\varphi V =- \alpha \xi + \beta \varphi V.
          \end{array}
  \right.
\end{equation}
Using formula \ref{DifTrans-Sas}, we get
$$ \alpha = \sigma= \varphi V(\rho) =0,\qquad \beta={\rm e}^{\rho}\quad and \quad {\rm div} V =2{\rm e}^{\rho}.$$
Which completes the proof.
\end{proof}
From this theorem, we can extract the following case:
\begin{corollary}
If $\rho=0$ then,  $(M,\overline{\varphi}, \psi,  \omega, g)$  is a Kenmotsu manifold if and only if ${\rm div} \psi =2$.
\end{corollary}
\begin{remark}
In this corollary, the result  was gotten  in \cite{BBB}.
\end{remark}

For our second application, let's define a $(1,1)$-tensor field $\tilde{\varphi}$ on $M$ by
$$ \widehat{\varphi}X=\varphi V \wedge X ,$$
that is, for all $X$ and $Y$  vector fields on $M$
\begin{eqnarray*}
\widehat{\varphi}X &=&g( \varphi^2 V, X)\xi - \eta(X)\varphi^2 V\\
&=&\eta(X) V - \theta^1(X)\xi.
\end{eqnarray*}
\begin{proposition}\label{ACMS2}
Let $(M, \varphi, \xi,  \eta, g)$ be a three dimensional  corner manifold. Then,  $(M,\widehat{\varphi}, \varphi V,  \theta^2, g)$
is an almost contact metric manifold.
\end{proposition}
\begin{proof}
One can adapt Proposition \ref{ACMS1}.
\end{proof}
\begin{thm}\label{ThCos}
The manifold $(M,\widehat{\varphi}, \varphi V,  \theta^2, g)$ is a cosymplectic manifold if and only if ${\rm div} V ={\rm e}^{\rho}$ and $\sigma = \varphi V(\rho) =0$.
\end{thm}
\begin{proof}
From Proposition \ref{PropِConnex}, we have
\begin{equation}\label{SysSas1}
 \left\{
          \begin{array}{lll}
 \nabla_\xi \varphi V = -\sigma V,\\
 \nabla_V  \varphi V =-\varphi V(\rho) V,\\
 \nabla_{\varphi V}\varphi V =\big( {\rm e}^{\rho} - {\rm div} V \big) V.
          \end{array}
  \right.
\end{equation}
On the other side, we have
\begin{equation}\label{SysSas2}
 \left\{
          \begin{array}{lll}
 -\alpha \widehat{\varphi} \xi - \beta \widehat{\varphi}^2\xi= \beta \xi -  \alpha V,\\
 -\alpha \widehat{\varphi}V - \beta \widehat{\varphi}^2 V = \alpha \xi + \beta V,\\
 -\alpha \widehat{\varphi}\varphi V - \beta \widehat{\varphi}^2 \varphi V =0.
          \end{array}
  \right.
\end{equation}
Using formula \ref{DifTrans-Sas}, we get
$$ \alpha =  \sigma = \beta=\varphi V(\rho) =0 \quad and \quad {\rm div} V ={\rm e}^{\rho}.$$
That is, $(M,\widehat{\varphi}, \varphi V,  \theta^2, g)$ is a cosymplectic manifold if and only if ${\rm div} V ={\rm e}^{\rho}$ and $\sigma=\varphi V(\rho) =0$.
\end{proof}
\begin{remark}
In this Theorem, the result  for $\rho=0$ was gotten  in \cite{BBB}.
\end{remark}
\begin{example}
In the class of examples given in section 4,  taking all the conditions into consideration with $\kappa = \kappa(x_2, x_3)$, $ \mu = \mu(x_2, x_3)$ and $\tau=\tau(x_1 , x_2)$, we get a  corner structure  $(\varphi ,\xi ,\eta, g)$ with 
 $$ \omega = \frac{\tau_2}{\tau} dx_2,\qquad  \psi =  \frac{\tau_2}{\tau \kappa^2 }\frac{\partial}{\partial x_2} \qquad and \qquad{\rm e}^{\rho}=\vert \psi \vert =   \frac{\tau_2}{\tau \kappa}.$$ 
In which case, one can get
$$ \sigma = g(\nabla_{\xi}V, \varphi V) =0 \qquad and \qquad \varphi V(\rho)= \frac{-2 \tau_2^2 \kappa_3}{\tau^2 \kappa^3 \mu}$$ 
 where $V= {\rm e}^{-\rho} \psi = \frac{1}{\kappa}\frac{\partial}{\partial x_2}=e_2$.
 So that $ \sigma = \varphi V(\rho)=0$  it suffices to take $\kappa_3 =0$ i.e. $\kappa = \kappa(x_2)$. Then, we get
 $$ \rm div V = \frac{(\tau \mu)_2}{\tau \kappa \mu}.$$
Thus, we can observe our cases as follows:\\
1)\; If $\frac{\tau_2}{\tau}=\frac{\mu_2}{\mu}$ we get $ {\rm div}V=2{\rm e}^{\rho}$ then, according to the Theorem  \ref{ThKen},   $(M,\overline{\varphi}, V,  \theta^1, g)$ is a $\beta$-Kenmotsu manifold with $\beta =\frac{\tau_2}{\tau \kappa}$.\\
 2)\; If $ \mu_2 =0$ we get $ {\rm div}V={\rm e}^{\rho}$ then, according to the Theorem \ref{ThCos},   $(M,\widehat{\varphi}, \varphi V,  \theta^2, g)$ is a cosymplectic manifold.

\end{example}
\subsection{From corner structure to normal almost contact metric structures}
The construction of normal almost contact metric structures (Sasakian, Kenmotsu, cosymplectic) from a corner structures which are not normal on a given 3-dimensional manifold M, in general, a non-trivial problem. The results obtained in the previous paragraph are a beautiful coincidence. Unfortunately, the techniques used do not enable us to obtain a Sasakian structure at all.

In 1992, J. C. Marrero \cite{MAR} proved that with certain deformation, we can get a trans-Sasakian structure starting from a Sasakian one. In \cite{AC}, generalized D-conformal deformations are applied to trans-Sasakian manifolds where the covariant derivatives of the deformed metric is evaluated under the condition that the functions used in deformation depend only on the direction of the characteristic vector field of the trans-Sasakian structure.  Other similar deformations are studied in \cite{BB, BM}.

Here, we present a certain deformation on the corner manifold which touches the structural tensor $\varphi $ and the metric tensor $g$ at the same time, which allow us to define new relations between all 3-dimensional almost  contact metric structures.

Let $( \varphi, \xi, \eta, g)$ be a corner structure on $M^{3}$. For any $X,$ $Y$ vector fields on $M$,  we mean a change of structure tensors of the form
\begin{equation}\label{NewStructure}
\left\{
	\begin{array}{llll}
		\tilde{\varphi}X = \varphi X + \theta^1( X)\xi,\\
		\tilde{\xi}=\xi,\\
		\tilde{\eta}=\eta - \theta^2,\\
		\tilde{g}(X, Y) = f g(X, Y) - f \eta(X)\eta(Y) + \tilde{\eta}(X)\tilde{\eta}(Y),
	\end{array}
	\right.
\end{equation}
where  $f$ a positive function  on $M$.
\begin{proposition} 
	The structure $(\tilde{\varphi} , \tilde{\xi} , \tilde{\eta} , \tilde{g})$  is an almost contact metric structure.
\end{proposition}
\begin{proof}
	The proof follows by a usual calculation, by using (\ref{CondPresqCont}).
\end{proof}

We denote the tensor field of type $(1,2)$ by $\tilde{N}^{(1)}$ on $M$ defined for any $X,$ $Y$ on $M$ by
$$
	\tilde{N}^{(1)}(X,Y)=[\tilde{\varphi},\tilde{\varphi}](X,Y)+2{\rm d}\tilde{\eta}(X,Y)\xi,
$$
where
$$ [\tilde{\varphi},\tilde{\varphi}](X,Y)=\tilde{\varphi}^{2}[X,Y]+[\tilde{\varphi} X,\tilde{\varphi} Y]-\tilde{\varphi}[\tilde{\varphi} X,Y]-\tilde{\varphi}[X,\tilde{\varphi} Y] .$$
By long direct calculation,  using (\ref{NewStructure}) one can get
\begin{eqnarray}\label{Ntilde0}
	\tilde{N}^{(1)}(X,Y)&=& N^{(1)}(X,Y) + \theta^2 (N^{(1)}(X,Y))\xi  \nonumber\\
	&-& \theta^2(\varphi X) \Big( N^{(3)}(Y) + \theta^2 \big( N^{(3)}(Y) \big)\xi \Big) 
	+ \theta^2(\varphi Y) \Big( N^{(3)}(X) + \theta^2 \big( N^{(3)}(X) \big)\xi \Big) \nonumber\\
	&-& 2 {\rm d}\theta^2(X,Y)\xi +  2 {\rm d}\theta^2 \big( \tilde{\varphi}X,\tilde{\varphi}Y \big) \xi 
\end{eqnarray}
with  $N^{(3)}$ is a tensor field on $M$ given  by
$$N^{(3)}(X)= \big( L_{\xi}\varphi )(X) = \varphi [X , \xi] - [\varphi X , \xi],$$
where $ L_{\xi}$ denotes the Lie derivative with respect to the vector field $\xi$.
\begin{proposition} 
	Let $( \varphi, \xi, \eta, g)$ be a corner structure on $M$. 
	The almost contact metric structure $(\tilde{\varphi} , \tilde{\xi} , \tilde{\eta} , \tilde{g})$  is normal if and only if
         $ \sigma = {\rm e}^{\rho}$.
\end{proposition}
\begin{proof}
Suppose that $( \varphi, \xi, \eta, g)$ is a corner structure on $M$. Then, we have $N_{\varphi}=0$, which gives
\begin{align*}
0 &= N_{\varphi}(X, \xi)= \varphi^2 [X , \xi] - \varphi[\varphi X , \xi],
\end{align*}
applying $\varphi$ we get
\begin{align*}
-\varphi [X , \xi] + [\varphi X , \xi] - \eta\big([\varphi X , \xi] \big)\xi =0,
\end{align*}
i.e. 
\begin{align*}
N^{(3)}(X) = -2{\rm d} \eta(\varphi X , \xi)\xi.
\end{align*}
So, (\ref{Ntilde0}) becomes
\begin{align}\label{Ntilde1}
	\tilde{N}^{(1)}(X,Y)=& 2{\rm d} \eta(X,Y)\xi 
	+ 2 \theta^2(\varphi X)  {\rm d} \eta(\varphi Y , \xi)\xi 
	-2 \theta^2(\varphi Y) {\rm d} \eta(\varphi X , \xi)\xi  - 2 {\rm d}\theta^2(X,Y)\xi\nonumber\\
	&+ 2 {\rm d}\theta^2 ( \varphi X, \varphi Y ) \xi +  2 \theta^2(\varphi X) {\rm d}\theta^2 ( \xi, \varphi Y ) \xi 
	+2 \theta^2(\varphi Y) {\rm d}\theta^2 ( \varphi X , \xi ) \xi . 	
\end{align}
From Proposition \ref{Propِd1Forms}, we have
\begin{align}\label{dtheta2}
{\rm d} \theta^2=  \sigma {\rm e}^{- \rho} {\rm d}\eta  +\frac{1}{2}({\rm e}^{\rho} -  {\rm div}V) \phi. 
\end{align}
By substituting (\ref{dtheta2}) in (\ref{Ntilde1}), we obtain
\begin{align}\label{Ntilde2}
\tilde{N}^{(1)}(X,Y)=&  2(1-\sigma {\rm e}^{-\rho}) \big( {\rm d} \eta(X,Y) -{\rm d} \eta(\varphi X , \xi) \theta^2(\varphi Y)
	-{\rm d} \eta(\xi, \varphi Y)\theta^2(\varphi X) \big) \xi. 	
\end{align}
Assume that
\begin{align}\label{eqNtilde}
{\rm d} \eta(X,Y) -{\rm d} \eta(\varphi X , \xi) \theta^2(\varphi Y)
	-{\rm d} \eta(\xi, \varphi Y)\theta^2(\varphi X)=0.
	\end{align}
	Putting $Y=\xi$, we get $ {\rm d} \eta(X,\xi)=0$ for all $X$ vector field on $M$. Then, the equation (\ref{eqNtilde}) implies 
	${\rm d} \eta =0$, which is contradiction. Hence, $\sigma= {\rm e}^{\rho}$.
This completes the proof.
\end{proof}
Now, the fundamental $2$-form  $\tilde{\phi}$ of $(\tilde{\varphi} , \tilde{\xi} , \tilde{\eta} , \tilde{g})$ is 
$$\tilde{\phi}( X,Y)= \tilde{g}(X, \tilde{\varphi} Y),$$

one can easily obtain
\begin{equation}\label{Omegah}
	\tilde{\phi} =  f \phi
\end{equation}

and hence, with the help of (\ref{dtheta2}), we have
\begin{equation}\label{Sys3}
	\left\{
	\begin{array}{lll}
		\tilde{\eta}=\eta- \theta^2\\
		\tilde{\phi} =  f \phi  
		
	\end{array}
	\right.
	\qquad \Rightarrow \qquad
	\left\{
	\begin{array}{lll}
		{ \rm d}\tilde{\eta}=(1 - \sigma {\rm e}^{- \rho}) {\rm d}\eta  +\frac{1}{2f}({\rm e}^{\rho} -  {\rm div}V)\tilde{\phi} \\
		{ \rm d}\tilde{\phi} =  { \rm d} (\ln f ) \wedge \tilde{\phi}.
		
	\end{array}
	\right.
\end{equation}
\begin{lemma}\label{Lem1}
	For any $3$-dimensional almost contact metric  manifold $(M,\tilde{\varphi} , \tilde{\xi} , \tilde{\eta} , \tilde{g})$, we have
	\begin{equation}
		{ \rm d} (\ln f ) \wedge \tilde{\phi} = \xi(\ln f ) \tilde{\eta} \wedge \tilde{\phi}.
	\end{equation}
\end{lemma}
\begin{proof}
	Let $\{\tilde{e}_0=\xi, \tilde{e}_1, \tilde{e}_2\}$  be the frame of vector fields  and $\{\tilde{\theta}^0=\tilde{\eta}, \tilde{\theta}^1, \tilde{\theta}^2\}$ be the dual frame of differential $1$-forms  on $M$. Then, 
	$$ \tilde{\phi}= 2 \tilde{e}_2 \wedge \tilde{e}_1,$$
	and 
	$$ { \rm d} (\ln f ) = \xi(\ln f)\tilde{\eta} + \tilde{\theta}^1(\ln f) \tilde{e}_1 + \tilde{\theta}^2(\ln f) \tilde{e}_2.$$
	Thus
	$${ \rm d} (\ln f ) \wedge \tilde{\phi} = \xi(\ln f ) \tilde{\eta} \wedge \tilde{\phi}.$$
\end{proof}

From (\ref{Sys3}) and Lemma \ref{Lem1}, we get
\begin{equation}\label{Sys4}
		\left\{
	\begin{array}{lll}
		{ \rm d}\tilde{\eta}=(1 - \sigma {\rm e}^{- \rho}) {\rm d}\eta  +\frac{1}{2f}({\rm e}^{\rho} -  {\rm div}V) \tilde{\phi} \\
		{ \rm d}\tilde{\phi} = \xi (\ln f ) \tilde{\eta} \wedge \tilde{\phi}.
		
	\end{array}
	\right.
\end{equation}
Based on these facts, we give the following Theorem:
\begin{thm}\label{ThSas}
Let $( M, \varphi, \xi, \eta, g)$ be a 3-dimensional corner manifold. The manifold  $(M, \tilde{\varphi} , \tilde{\xi} , \tilde{\eta} , \tilde{g})$ define above is a Trans-Sasakian manifold of type
$$\Big(\tilde{\alpha}=\frac{1}{2f}({\rm e}^{\rho} -  {\rm div}V) , \tilde{\beta}=\frac{1}{2}\xi (\ln f ) \Big)$$
if and only if $ \sigma = {\rm e}^{\rho}$.
\end{thm}
Here, we mention the three interesting cases as follows:
\begin{corollary}
Let $( M, \varphi, \xi, \eta, g)$ be a 3-dimensional corner manifold with $\sigma = {\rm e}^{\rho}$. Then, the manifold  $(M, \tilde{\varphi} , \tilde{\xi} , \tilde{\eta} , \tilde{g})$ define above is:
\begin{equation}\label{DifStructures1}
	\left\{
	\begin{array}{lll}
		(a): \;  Sasakian \quad iff \quad {\rm div}V = {\rm e}^{\rho}-2f \quad and \quad \xi(f)=0.\\
		(b): \;  Kenmotsu \quad iff \quad {\rm div}V = {\rm e}^{\rho}\quad and \quad \xi(f)=2f.\\
		(c): \;  Cosymplectic \quad iff \quad {\rm div}V = {\rm e}^{\rho}\quad and \quad \xi(f)=0.\\
	\end{array}
	\right.
\end{equation}
\end{corollary}
Consequently, for $\rho=0$ we get the following:
\begin{corollary}
Let $( M, \varphi, \xi, \eta, g)$ be a 3-dimensional unit corner manifold with $\sigma = 1$ Then, the manifold  $(M, \tilde{\varphi} , \tilde{\xi} , \tilde{\eta} , \tilde{g})$ define above is:
\begin{equation}\label{DifStructures2}
	\left\{
	\begin{array}{lll}
		(a): \;  Sasakian \quad iff \quad {\rm div}\psi = 1-2f \quad and \quad \xi(f)=0.\\
		(b): \;  Kenmotsu \quad iff \quad {\rm div}\psi = 1\quad and \quad \xi(f)=2f.\\
		(c): \;  Cosymplectic \quad iff \quad {\rm div}\psi = 1\quad and \quad \xi(f)=0.\\
	\end{array}
	\right.
\end{equation}
\end{corollary}
\begin{proposition}
The results of the Theorem \ref{ThSas} do not hold if $\omega$ is closed.
\end{proposition}
\begin{proof}
We have
\begin{align*}
\sigma &= g(\nabla_{\xi}V, \varphi V)\\
&= {\rm e}^{-2 \rho}g(\nabla_{\xi}\psi, \varphi \psi),
\end{align*}
if $\omega$ is closed, that is $ g(\nabla_{\xi}\psi, \varphi \psi) = g(\nabla_{\varphi \psi}\psi, \xi)$. Then,
\begin{align*}
\sigma &=- {\rm e}^{-2 \rho}g(\psi, \nabla_{\varphi \psi}\xi)\\
&=0.
\end{align*}
But, in Theorem \ref{ThSas}, $\sigma = {\rm e}^{\rho} \neq0$.
\end{proof}
\textbf{Open question:}\\
The last proposition poses an open question on how to obtain a Sasakian manifold from a corner manifold with a closed differential form $\omega$?.

\end{document}